\documentclass[draft]{article}

\usepackage{amsmath}
\usepackage{amsthm}
\usepackage{amscd}
\usepackage{amssymb}
\usepackage{latexsym}

\newcommand{\LL}{{\cal L}}
\newcommand{\CA}{{\cal A}}
\newcommand{\CB}{{\cal B}}

\newcommand{\LG}{\mathfrak L}

\newcommand{\vbar}{|}
\newcommand{\rvbar}{|}
\newcommand{\chop}{\dag}

\def\epsilon{\varepsilon}

\newcommand{\Out}{\mbox{Out}}
\newcommand{\Aut}{\mbox{Aut}}
\newcommand{\FN}{F_{N}}
\newcommand{\inv}{^{-1}}
\newcommand{\BBT}{\mbox{BBT}}

\newcommand{\R}{\mathbb R}
\newcommand{\Z}{\mathbb Z}

\newcommand{\N}{\mathbb N}
\newcommand{\Hy}{\mathbb H}
\newcommand{\A}{{\cal A}}

\newcommand{\FM}{F_{M}}
\def\phi{\varphi}

\def\qed{\hfill\rlap{$\sqcup$}$\sqcap$\par}

\def\tilde{\widetilde}

\newtheorem{thm}{Theorem}[section]

\newtheorem{lem}[thm]{Lemma}
\newtheorem{prop}[thm]{Proposition}

\theoremstyle{definition}
\newtheorem{defn}[thm]{Definition}
\newtheorem{rem}[thm]{Remark}
\newtheorem{quest}[thm]{Question}

\newtheorem*{acknowledgements}{Acknowledgements}

\newcommand{\affiliationone}{}
\newcommand{\email}{\\}
\renewcommand{\and}{ and }

\title{$\R$-trees and laminations for free groups I: Algebraic laminations}

\author{Thierry Coulbois, Arnaud Hilion \and Martin Lustig}

\begin{document}

\maketitle

\section{Introduction}

This paper is the first of a sequence of three papers, where the
concept of an $\R$-tree dual to (the lift to the universal covering of)
a measured geodesic lamination $\LG$ in a hyperbolic surface $S$ is
generalized to arbitrary $\R$-trees provided with a
(very small) action of the free group $\FN$ of finite
rank $N \geq 2$ by isometries.

\smallskip

In \cite{chl1-II} to any such $\R$-tree
$T$ a {\em dual algebraic lamination} $L^{2}(T)$ is associated in
a meaningful way, and in \cite{chl1-III} we consider
invariant measures (called {\em currents}) $\mu$ on $L^{2}(T)$ and
investigate the
induced {\em dual metric} $d_{\mu}$ on $T$.

\smallskip

In this first paper we define and study the basic tools for the two
subsequent papers: {\em laminations} in the free group $\FN$.  We will
use three different approaches, {\em algebraic laminations} $L^{2}$,
{\em symbolic laminations} $L_{\CA}$, and {\em laminary languages}
$\LL$. Each of them will be explained in detail, and each has its own
virtues.  Algebraic laminations do not need a specified basis of $\FN$
and are hence of conceptional superiority. The other two objects are
concretely defined in terms of infinite words (for symbolic
laminations) or of finite words (for laminary languages) in a fixed
basis $\CA$. They are more practical for many tasks: Symbolic
laminations are more suited for dynamical and laminary languages more
for combinatorial purposes.  The set of each of these three objects
come naturally with a topology, a partial order, and an action by
homeomorphisms of the group $\Out(\FN)$ of outer automorphisms of
$\FN$.  We will prove that the three approaches are equivalent:

\begin{thm}
\label{theoremone}
Let $\FN$ denote the free group of finite rank $N \geq 2$, and let
$\CA$ be a basis of $\FN$.  There are canonical
$\Out(\FN)$-equivariant, order preserving homeomorphisms
\[
\Lambda^{2}(\FN) \longleftrightarrow \Lambda_{\CA}
\longleftrightarrow
\Lambda_{\LL}(\CA)
\] 
between the space $\Lambda^{2}(\FN)$ of algebraic laminations in
$\FN$, the space $\Lambda_{\CA}$ of symbolic laminations in $\CA^{\pm
1}$, and the space $\Lambda_{\LL}(\CA)$ of laminary languages in
$\CA^{\pm 1}$.
\end{thm}

Symbolic laminations are subshifts (= symbolic flows) as classically
used in symbolic dynamics, except that we work with the free group
$\FN = F(\CA)$ rather than with the free monoid $\CA^{*}$.  Similarly,
laminary languages over the alphabet $\cal A$ rather than $\CA^{\pm
1}=\CA \cup \CA^{-1}$ are already studied in combinatorics, compare
for instance \cite{narbel}.

\smallskip

As in the surface case, the subset $\Lambda_{\mbox{\scriptsize rat}}
\subset \Lambda^2(\FN)$ of rational laminations, each corresponding to
a finite collection of non-trivial conjugacy classes in $F_{N}$ (see
\S\ref{sec:alglam}), is of special interest.  Contrary to the
analogous statement for measured laminations on a surface, or for
currents on $\FN$ (compare \cite{mart}), we obtain in the setting of
algebraic laminations:

\begin{thm}\label{thm:ratdense}
Rational laminations are not dense in $\Lambda^2(\FN)$. However,
the closure $\overline{\Lambda}_{\mbox{\scriptsize rat}}$ contains
all minimal laminations.
\end{thm}

Algebraic laminations, as defined and studied in this paper, have
three direct ``ancesters", all three of them inspired by geodesic
laminations on surfaces: In \cite{luhabil} {\it combinatorial
laminations} are defined to study decomposable automorphisms of
$F_{N}$, in \cite{bfhtits1} an {\it attracting lamination} is
associated to each exponential stratum of an automorphism of $\FN$
(see \S\ref{sec:alglam}), and in \cite{ll4} a kind of laminations is
associated to certain $\R$-tree actions of $\FN$.

\smallskip

This paper (as well as the subsequent ones \cite{chl1-II} and
\cite{chl1-III}) is a further attempt to bridge the ``cultural gap''
between two mathematical communities: symbolic and combinatorial
dynamics on one hand, and geometric group theory on the other.  Notice
that in geometric group theory the notion of an algebraic lamination
extends naturally to the more general setting of word-hyperbolic
groups.

\smallskip

We hope to have given enough detail to carry along the novice reader
from the ``other'' mathematical subculture, and not too much to bore
the expert reader from ``this'' one.

\bigskip

\begin{acknowledgements}
This paper originates from a workshop organized at the CIRM in April
05, and it has greatly benefited from the discussions started there
and continued around the weekly Marseille seminar ``Teichm\"uller''
(partially supported by the FRUMAM).
\end{acknowledgements}

%%%%%%%%%%%%%%%%%%%%%%%%%%%%%%%%%%%%%%%%%%%%%%%%%%%%%%%%%%%%

\section{Algebraic laminations}\label{sec:alglam}

Let $\FN$ denote the free group of finite rank $N \geq 2$, and let
$\partial \FN$ denote its Gromov boundary, as usual equipped with the
action of $\FN$ (from the left) and with Gromov's topology at
infinity, which gives $\partial \FN$ the topology of a Cantor set.
The choice of a basis $\cal A$ of $\FN$ allows us to identify the
elements of $\FN$ with reduced words $w = x_{1} x_{2} \ldots x_{n}$
(with $x_{i+1} \neq x_{i}\inv$) in ${\cal A}\cup{\cal A}^{-1}$, and
thus defines in particular the length function $w \mapsto |w|_{\CA} =
n$ on $\FN$.  This length function induces the {\em word metric}
$d_{\CA}(v, w) = | v\inv w |_{\CA}$ on $\FN$, which in turn defines a
metric on $\partial \FN = \{ x_{1} x_{2} x_{3} \ldots \mid x_{i} \in
\CA^{\pm 1}, x_{i+1} \neq x_{i}\inv \}$, stated explicitely in
\S\ref{subsec:topologyonLambda}.

Choosing another basis gives rise to a Lipschitz-equivalent metric on
$\FN$ and to a H\"older-equivalent metric on $\partial \FN$ (compare
\cite{gh}).  As a consequence, the topology on $\FN \cup \partial \FN$
induced by the word metric does not depend on the choice of the basis
$\cal A$.  More details are given below in
\S\ref{sec:outfnactiononLambda}.  Note that $\FN \cup \partial \FN$ as
well as $\partial \FN$ are compact spaces, and that every $\FN$-orbit
in $\partial \FN$ is dense.

\smallskip

For any element $w \neq 1$ of $\FN$ we denote by $w^{+\infty}$ the
limit in $\partial \FN$ of the sequence $(w^n)_{n\in\N}$ and by
$w^{-\infty}$ that of $(w^{-n})_{n\in\N}$.  If $w = x_{1} \ldots x_{p}
\cdot y_{1} \ldots y_{q} \cdot x_{p}\inv \ldots x_{1}\inv$ is a
reduced word in $\CA^{\pm 1}$, with $y_{q} \neq y_{1}\inv$, then
\[
w^{+\infty} = x_{1} \ldots x_{p}
\cdot y_{1} \ldots y_{q} \cdot y_{1} \ldots y_{q}
\cdot  y_{1} \ldots y_{q} \cdot \ldots
\]

\smallskip

Following standard notation (see for example \cite{kapo1,kapo2}), we
define
\[
\partial^{2}\FN = \partial \FN \times \partial \FN \smallsetminus
\Delta
\, ,
\]
where $\Delta$ denotes the diagonal in $\partial \FN \times \partial
\FN$. It follows directly that $\partial^{2} \FN$ inherits from
$\partial \FN$ a topology and an $\FN$-action, given by $w(X, Y) =
(wX, wY)$ for any $w \in \FN$ and any $X,Y \in \partial \FN$ with $ X
\neq Y$. The set $\partial^{2} \FN$ admits also the {\em flip}
involution $(X, Y) \mapsto (Y, X)$, which is an $\FN$-equivariant
homeomorphism.  Note that $\partial^{2}\FN$ is not compact.

\smallskip

\begin{defn}\label{def:alglam} 
An {\em algebraic lamination} is a subset $L^{2}$ of $\partial^{2}\FN$
which is non-empty, closed, symmetric (= flip invariant) and
$\FN$-invariant. The set of all algebraic laminations is denoted by
$\Lambda^2 = \Lambda^2(\FN)$.
\end{defn}

The set $\Lambda^2$ of algebraic laminations inherits naturally a
Hausdorff topology from $\partial^{2} \FN$ which we will discuss in
\S\ref{subsec:topologyonLambda}.

\smallskip

In \cite{bfhtits1}, M.~Bestvina, M.~Feighn and M.~Handel associate an
attracting lamination to each exponential stratum of an automorphism
of $\FN$. These laminations are laminations in our sense.  However, in
\cite{bfhtits1} there is no topology introduced on the space of
laminations but rather only on $\partial^2\FN$, and even there, their
topology differs slighty from ours.

\smallskip

An important special class of algebraic laminations are the {\em
rational} laminations, which are finite unions of {\em minimal
rational} laminations $L(w)$, defined for any $w \in \FN
\smallsetminus \{1\}$ by:
\[
L(w) = \{(vw^{-\infty}, vw^{+\infty}) \, \mid \, v \in
\FN \}
\cup \,
\{(vw^{+\infty}, vw^{-\infty}) \, \mid \, v \in
\FN \}
\]
Note that the lamination $L(w)$ depends only on the conjugacy class of
$w$.  We denote by $\Lambda_{\mbox{\scriptsize rat}}$ the subspace of
rational laminations.  The Hausdorff topology on $\Lambda^{2}$ is
stronger than one might intuitively expect.  In particular on obtains
the following result, proved in \S\ref{subsec:topologyonLambda}:

\begin{prop}\label{prop:notdense}
The subset $\Lambda_{\mbox{\scriptsize rat}}$ is not dense in
$\Lambda^{2}$.
\end{prop}

We observe that there is a natural (left) action of $\Out(\FN)$ on
$\Lambda^{2}$, induced by the action of $\Aut(\FN)$ on $\partial
\FN$. Indeed, an automorphism of $\FN$ is a bi-Lipschitz homeomorphism
on $\FN$ and extends continuously to the boundary. Inner automorphisms
act by left-multiplication on the boundary and thus trivially on the
space $\Lambda^{2}$ of algebraic laminations (as the latter are
$\FN$-invariant subsets of $\partial^2 \FN$).  More details about the
$\Out(\FN)$-action on $\Lambda^{2}$ will be given in
\S\ref{sec:outfnactiononLambda}.

Note that this action restricts to an action of $\Out(\FN)$ on the
space of rational laminations $\Lambda_{\mbox{\scriptsize rat}}$: If
$\alpha$ is an automorphism of $\FN$ and $\widehat\alpha$ its class in
the outer automorphism group $\Out(\FN)$ and, if $w$ is an element of
$\FN$, $\widehat\alpha(L(w))=L(\alpha(w))$.

\smallskip

To stimulate the interest of the reader in these rather delicate
matters we would like to pose here a question which is inspired by the
thesis of R.~Martin \cite{mart}:

\begin{quest}
\label{uniqueminimal}
Let $\CA$ be any basis of $\FN$, and fix $a\in {\cal A}$ arbitrarily.
Is the closure $\overline{\Out(\FN) L(a)}$ of the $\Out(\FN)$-orbit of
$L(a)$ a minimal closed $\Out(\FN)$-invariant non-empty subset of
$\Lambda^2$ ?  If so, is it the unique such {\em minimal} set?
\end{quest}

An answer to this question will be given in Proposition \ref{thierry}.
Note that if $N=2$ and $\{ a,b\}$ is a basis of $F_2$ and $[a,b]=a\inv
b\inv ab$, then it is well known that for any automorphism $\alpha$ of
$\FN$, $\alpha([a,b])$ is conjugated to either $[a,b]$ or its
inverse. Therefore $L([a,b])$ is a global fixed point of the action of
$\Out(\FN)$ on $\Lambda$.

%%%%%%%%%%%%%%%%%%%%%%%%%%%%%%%%%%%%%%%%%%%%%%%%%%%

\section{Surface laminations}\label{section:surfacelamination}

An important class of algebraic laminations comes from geodesic
laminations on hyperbolic surfaces.  The discussion started below, to
compare algebraic laminations in general with laminations on surfaces,
is carried further in \cite{chl1-II} and \cite{chl1-III}.  Throughout
this section we assume a certain familiarity of the reader with this
subject; for background see for example \cite{cassonbleiler} and
\cite{flp}. Note that this section can be skipped by the reader
without loss on the intrinsic logics of the material presented in this
paper.

\smallskip

Let $S$ be a hyperbolic surface with non-empty boundary and negative
Euler characteristic, and fix an identification $\pi_{1} S = \FN$.
The surface $S$ is provided with a hyperbolic structure, given by an
identification of the universal covering $\tilde S$ with a convex part
of the hyperbolic plane $\Hy^{2}$, which realizes the deck
transformation action of $F_{N} = \pi_{1}S$ on $\tilde S$ by
hyperbolic isometries.  Let ${\LG}$ be a geodesic lamination on $S$
and let $\tilde {\LG}$ be the (full) lift of $\LG$ to the universal
covering $\tilde S$ of $S$.  The induced identification (an
$\FN$-equivariant homeomorphism!)  between $\partial \FN$ and the
boundary at infinity $\partial \tilde S$ of $\tilde S$ defines for any
leaf $l$ of $\tilde {\LG}$ a pair of endpoints $(X, Y) \in
\partial^{2}\FN$, as well as its flipped pair $(Y, X)$. The set of all
such pairs is easily seen to define (via the above identification
$\partial \FN = \partial \tilde S$) an algebraic lamination
$L^{2}({\LG}) \in\Lambda^2(\FN)$.

\begin{defn}
An algebraic lamination $L^{2} \in \Lambda^{2}(\FN)$ is called an {\em
algebraic surface lamination} if there exists a hyperbolic surface $S$
and an identification $\pi_{1} S = \FN$ such that for some geodesic
lamination $\LG$ on $S$ one has:
\[
L^{2} = L^{2}({\LG})
\]
\end{defn}

At first guess it may seem that the space $\Lambda^2(\FN)$ is a rather
weak analogue of the space of geodesic laminations in a surface.
Notice however that, if $L^{2} \in \Lambda^{2}(\FN)$ is an algebraic
surface lamination with respect to an isomorphism $\pi_{1} S_{1} =
\FN$ for some surface $S_{1}$, and if $S_{2}$ is a second surface with
identification $\pi_{1} S_{2} = \FN$, then typically a biinfinite
geodesic on $S_{2}$, which realises an element of $L^{2}$, will
self-intersect: Thus $L^{2}$ does not admit a realization as geodesic
lamination on $S_{2}$.

%%%%%%%%%%%%%%%%%%%%%%%%%%%%%%%%%%%%%%%%%%%%%%%%%%%%%%%%%%%

\section{Symbolic laminations}

To a basis $\cal A$ there is naturally associated the space
$\Sigma_{\cal A}$ of biinfinite reduced words $Z$ in ${\cal A} \cup
{\cal A}^{-1}$ with letters indexed by $\Z$:
\[
\Sigma_{\cal A} = \{ Z = \ldots z_{i-1} z_{i} z_{i+1} \ldots \mid
z_{i} \in {\cal A} \cup {\cal A}^{-1}, z_{i} \neq z_{i+1}^{-1} \,\,
\hbox{\rm for all} \,\, i\in \Z \}.
\]
We want to stress that in this paper a biinfinite word comes always
with a $\Z$-indexing, i.e.  formally speaking, a biinfinite word is a
map $Z: \Z \to \CA \cup \CA\inv$.  For example, the non-indexed
``biinfinite word''
\[
\ldots a b a b a b \ldots
\]
becomes a biinfinte word $Z$ only after specifying $z_{1} = a$ or
$z_{1} = b$, which we indicate notationally by writing $Z = \ldots b a
b\cdot a b a \ldots$ or $Z = \ldots a b a \cdot b a b \ldots$
respectively.

\smallskip

As usual, $\Sigma_\CA$ comes with a canonical infinite cartesian
product topology that makes it a Cantor set, and with a shift operator
$\sigma: \Sigma_\CA \to \Sigma_\CA$, given by
\[
\sigma( Z) =  Z' \, ,
\]
where $ Z = \ldots z_{i-1} z_{i} z_{i+1} \ldots$ and $ Z' = \ldots
z'_{i-1} z'_{i} z'_{i+1} \ldots$ with $z'_{i} = z_{i+1}$.  Of course,
$\sigma$ is a homeomorphism.

\smallskip

For each biinfinite word $Z=\ldots z_{i-1} z_{i} z_{i+1}\ldots$ we
denote its {\em inverse} by 
\[
Z^{-1}=\ldots z'_{i-1} z'_{i}
z'_{i+1}\ldots\,,\mbox{ where }z'_i=(z_{1-i})^{-1}\, .
\]  
Again, the inversion
map $\Sigma_{\CA} \to \Sigma_{\CA}\, , \, \, Z \mapsto Z\inv$ is
easily seen to be a homeomorphism.  A subset $L$ of $\Sigma_\CA$ is
called {\em symmetric} if $L = L^{-1}$.

\begin{defn}\label{def:symblam}
A {\em symbolic lamination} in $\CA^{\pm 1}$ is a non-empty subset
$L_{\cal A} \subset \Sigma_{\cal A}$ which is closed, symmetric and
$\sigma$-invariant.  Together with the restriction of $\sigma$ to
$L_{\cal A}$ (which we continue to call $\sigma$) it is a {\em
symbolic flow}.  The elements of a symbolic lamination are sometimes
called the {\em leaves} of the lamination.  We denote the set of
symbolic laminations in $\CA^{\pm 1}$ by $\Lambda_{\CA}\,$.
\end{defn}

In symbolic dynamist's terminology, any symbolic lamination is a
subshift of the subshift of finite type on the alphabet $\CA \cup
\CA\inv$ which consists of all biinfinite reduced words.

As $\Sigma_\CA$ is compact and symbolic laminations are closed, we
get:

\begin{lem}\label{lem:decreasingintersection}
The intersection of a decreasing sequence
\[
L_{\CA}\supset L'_{\CA} \supset L''_{\CA} \supset \ldots 
\]
of symbolic laminations is a symbolic lamination. In particular it is
non-empty.\qed
\end{lem}

Once the basis ${\cal A}$ is fixed, every boundary point $X \in
\partial \FN$ corresponds canonically to a reduced, (one-sided)
infinite word $X = x_{1} x_{2} \ldots$ with letters in $\CA^{\pm 1}$.
For such a (one-sided) infinite word $X$ we denote by $X_n$ its prefix
(= initial subword) of length $n$.  For every pair $(X, Y) \in
\partial^{2} \FN$ we define a biinfinite reduced word
\[
X^{-1} Y = \ldots x_{k+2}^{-1} x_{k+1}^{-1} \cdot y_{k+1} y_{k+2}
\ldots\, ,
\]
where $X_{k} = x_{1} x_{2} \ldots x_{k} = y_{1} y_{2} \ldots y_{k} =
Y_{k}$ is the longest common prefix of $X$ and $Y$.

\smallskip

There is a subtlety in the last definition which we would like to
point out: Although for any $X \neq Y \in \partial \FN$ the biinfinite
(indexed) word $X\inv Y$ is well defined by our above definition, this
particular way to associate the indices from $\Z$ to the non-indexed
``biinfinite word'' $\ldots x_{k+2}^{-1} x_{k+1}^{-1} y_{k+1}
y_{k+2}\ldots$ is really in no way canonical, and often it does not
behave quite naturally, in particular with respect to the action of
$\Aut \FN$.  Indeed, a biinfinite symbol sequence, contrary to a
finite or a one-sided infinite one, doesn't really come by nature with
a canonical indexing, but rather corresponds to the whole
$\sigma$-orbit of a biinfinite word in $\Sigma_{\cal A}$.
Nevertheless one obtains as direct consequence of the definitions:

\begin{rem}
\label{noncontinuous}
The map $\rho_{\CA}:~\begin{array}[t]{rcl}\partial^{2}
\FN&\to&\Sigma_{\CA}\\ (X,Y)&\mapsto&X^{-1} Y
\end{array}$
is continuous.
\end{rem}

We note that the biinfinite indexed word from $\Sigma_\CA$ associated
via $\rho_{\CA}$ to $w(X, Y)$, for any $w \in \FN$, can differ from
the indexed word $X^{-1} Y$ only by an index shift.  Conversely, for
the pair $(X, Y) \in \partial^{2} \FN$ with maximal common initial
subword $X_{k} = Y_{k}$ as above, the map $\rho_{\CA}$ associates the
biinfinite indexed word $\sigma^m(X^{-1} Y)$ to the pair $Y_{k+m}^{-1}
(X, Y)$ for $m \geq 0$, and to $X_{k-m}^{-1} (X, Y)$ for $m \leq 0$.

\smallskip

Hence the map $\rho_{\CA}$ maps every $\FN$-orbit in $\partial^{2}
\FN$ onto a $\sigma$-orbit in $\Sigma_{\CA}$, and thus induces a well
defined map from $\FN$-orbits in $\partial^{2} \FN$ to $\sigma$-orbits
in $\Sigma_{\CA}$.  It is easy to see that this map between orbits is
bijective, and that, moreover, this bijection respects the topology on
both sides: Closed sets of $\FN$-orbits are mapped to closed sets of
$\sigma$-orbits, and conversely.  Finally, we note that the flip on
$\partial^{2} \FN$ corresponds to the inversion of biinfinite words in
$\Sigma_\CA$.

\smallskip

Thus, given $L^{2} \in \Lambda^{2}$, we can define a symbolic
lamination $L_{\cal A}$ by
\[
L_{\cal A} = \rho_{\CA}(L^{2}) = \{ X^{-1} Y \mid (X, Y) \in L^{2} \}.
\]

\smallskip

Conversely, given a symbolic lamination $L_{\cal A}$ as above, one
obtains an algebraic lamination $L^{2} = \rho_{\CA}^{-1}(L_{\cal A})$
which consists of all pairs $w(Z_{-}, Z_{+})$, for all $w \in \FN$,
and all $Z = \ldots z_{i-1} z_{i} z_{i+1} \ldots \in L_{\cal A}$ with
associated right-infinite words $Z_{-} = z_{0}^{-1} z_{-1}^{-1}
z_{-2}^{-1}\ldots$ and $Z_{+} = z_{1} z_{2} \ldots\, \, \,$.

We summarize the above discussion:

\begin{prop}
\label{prop:laminationflow}
For any basis $\cal A$ of the free group $\FN$, the maps $L^{2}
\mapsto L_{\cal A} = \rho_{\CA}(L^{2})$ and $L_{\cal A}\mapsto L^{2} =
\rho_{\CA}^{-1}(L_{\cal A})$ define a bijection
\[
\rho^{2}_{\CA}: \Lambda^2(\FN) \to \Lambda_{\CA}
\]
between the set $\Lambda^2(\FN)$ of algebraic laminations $L^{2}$ and
the set $\Lambda_{\CA}$ of symbolic laminations $L_{\CA}$ in $\CA^{\pm
1}$.  \qed
\end{prop}

The map $\rho^2_{\CA}$ respects the partial order given on algebraic
or symbolic laminations by the inclusion as subsets of
$\partial^{2}\FN$ or $\Sigma_{\CA}$ respectively. In particular, a
minimal lamination $L_{\CA}$ (or $L^{2}$) with respect to this partial
order is precisely given by the analogous property that characterizes
classically {\em minimal} symbolic flows: Every $< \sigma,
(\cdot)\inv>$-orbit (or $<\FN, \hbox{\rm flip}>$-orbit, respectively)
is dense in the lamination.  Moreover, we note that Lemma
\ref{lem:decreasingintersection} holds for algebraic laminations.

\smallskip

In order to connect the content (and also the notations) introduced in
this section to the already existing notions in symbolic dynamics, we
note:

A symbolic flow $\sigma: \Sigma_{0} \to \Sigma_{0}$ in the ``classical
sense", i.e. a symbolic flow only on the letters of $\cal A$ (and not
of ${\cal A}^{-1}$), gives directly rise to a symbolic lamination
$L_{\CA}(\Sigma_{0}) = \Sigma_{0} \cup {\Sigma_{0}}^{-1} \in
\Lambda_{\CA}$.  Conversely, a symbolic lamination $L_{\CA} \in
\Lambda_{\CA}$ or a symbolic flow $\sigma: L_{\CA} \to L_{\CA}$ is
called {\em orientable} if $L$ can be written as disjoint union
$L_{\CA} =L_{+} \cup L_{+}^{-1}$ of two $\sigma$-invariant closed
subsets $L_{+} $ and $L_{+}^{-1}$ that are inverses of each other, and
it is called {\em positive} if one of them, say $L_{+}$, only uses
letters from $\cal A$ (and not from ${\cal A}^{-1}$).

\begin{rem}
\label{lem:abhs}
The fact that the laminations considered are positive is crucial for
many of the traditional approaches and methods of symbolic
dynamics. Similarly, for laminations (or foliations) on surfaces,
almost always one first considers the orientable case and later tries
to pass to the general situation via branched coverings. Note that in
the context of free groups considered here any such attempt would miss
most of the typical phenomena, and that hence struggling with the
general kind of non-orientable laminations seems unavoidable. For an
interesting case of such an encounter of the free group environment
with the ``already existing culture'' in the context of the Rauzy
fractal see \cite{abhs}.
\end{rem}

%%%%%%%%%%%%%%%%%%%%%%%%%%%%%%%%%%%%%%%%%%%%%%%%%%%%%%%%%%%%

\section{Laminary languages}\label{subsec:lamlang}

As before, we fix a basis $\cal A$ of $\FN$, and we denote by $F({\cal
A})$ the set of reduced words in $\CA^{\pm 1}$. Although there is a
canonical identification between $\FN$ and $F({\cal A})$, it is
helpful in the context of this section to think of the elements of
$F({\cal A})$ as words and not as group elements.

\begin{defn}
\label{def:language}
Let $S$ be any (finite or infinite) set of finite, one-sided infinite
or biinfinite reduced words in $\CA^{\pm 1}$.  We denote by ${\cal
L}(S) \subset F({\cal A})$ the {\em language} generated by $S$, i.e.
the set of all finite subwords (= {\it factors}) of any element of
$S$.  Moreover, for any integer $n$ we denote by ${\cal L}_n(S)$ the
subset of ${\cal L}(S)$ consisting of words of length smaller or equal
to $n$.
\end{defn}

We specially have in mind the language associated to a (symbolic)
lamination. We thus abstractly define laminary languages which are in
one-to-one correspondence with (symbolic) laminations.

\begin{defn}
\label{def:lamlang}
A non-empty set ${\cal L} \subset F(\CA)$ of finite reduced words in
$\CA^{\pm 1}$ is a {\em laminary language} if it is (i) symmetric,
(ii) factorial and (iii) bi-extendable.  By this we mean that it is
closed with respect to (i) inversion, (ii) passing to subwords, and
(iii) that for any word $u \in \cal L$ there exists a word $v \in \cal
L$ in which $u$ occurs as subword other than as prefix or as suffix:
$v = w u w'$ is a reduced product, with nontrivial $w, w' \in F(\CA)$.
We denote by $\Lambda_{\LL} = \Lambda_{\LL}(\CA)$ the set of laminary
languages over a fixed basis $\CA$.
\end{defn}

It is obvious from the definition that the set $\Lambda_\LL$ is closed
under (possibly infinite) unions in $F(\CA)$, and also under nested
intersections (compare with Lemma \ref{lem:decreasingintersection}).
Note that the analogy of the former statement, for symbolic
laminations rather than laminary languages, is false: An infinite
union of symbolic laminations will in general not be a symbolic
lamination; one first needs to take again the closure in
$\Sigma_{\CA}$.  Note also that for any symbolic lamination $L_{\CA}
\subset \Sigma_{\CA}$ the language $\LL(L_{\CA})$ is laminary.

\smallskip

For an infinite language $\LL \subset F(\CA)$, we denote by $L({\cal
L})$ the set of all biinfinite words from $\Sigma_{\cal A}$ whose
finite subwords are subwords of elements from $\LL \cup \LL\inv$.  As
$\cal L$ is infinite (hence in particular, if $\LL$ is a laminary
language), the definition enforces that $L({\cal L})$ is not empty.
It follows directly that $L({\cal L})$ is indeed a symbolic
lamination.  We thus obtain a one-to-one correspondence between
symbolic laminations and laminary languages (always for a fixed basis
$\CA$ of $\FN$): For any symbolic lamination $L_{\CA}$ one has
\[
L({\cal L}(L_{\CA}))=L_{\CA} \, ,
\]
and conversely, for any laminary language ${\cal L}$ one has
\[
{\cal L}(L({\cal L}))={\cal L} \, .
\]
Moreover, a language ${\cal L}$ is laminary if and only if it is
infinite, and if the last equation holds.  For any set $S$ of finite,
one-sided infinite or biinfinite reduced words in $\CA^{\pm 1}$, where
we assume that $S$ is infinite in case $S \subset F(\CA)$, we observe
that ${\cal L}(L({\cal L}(S)))$ is the largest laminary language
contained in $\LL(S)$.  We call $L({\cal L}(S))$ the symbolic
lamination and ${\cal L}(L({\cal L}(S)))$ the laminary language {\em
generated} by $S$.  We summarize this discussion:

\begin{prop}
\label{bijectionsymbolic}
For any finite alphabet $\CA$ the maps $L_{\CA} \mapsto \LL(L_{\CA})$
and $\LL\mapsto L(\LL)$ define a bijection
\[
\rho^{\CA}_\LL:
\Lambda_{\CA}
\to \Lambda_{\LL}
\]
between the set $\Lambda_{\CA}$ of symbolic laminations $L_{\CA}$ and
the set $\Lambda_{\LL}$ of laminary languages $\LL$ in $\CA^{\pm 1}$.
\qed
\end{prop}

As in Proposition \ref{prop:laminationflow}, the bijection
$\rho^{\CA}_\LL: \Lambda_{\CA} \to \Lambda_{\LL}$ respects the partial
order given by the inclusion.

\smallskip

To enforce the link between symbolic laminations and their laminary
languages we introduce the following notation and state the following
lemma, which will be used in the sequel: For any integer $k \geq 0$
and any reduced word $w = x_{1} x_{2} \ldots x_{n}\in F(\CA)$ denote
by $w\chop_{k}$ (``w chop k'') the word

\begin{enumerate}
\item[(a)]
$w\chop_{k} = 1$, if $|w| \leq 2 k$, and 

\item[(b)]
$w\chop_{k} = x_{k+1} x_{k+2} \ldots x_{n-k}$, if $|w| > 2 k$.
\end{enumerate}

Similarly, for any integer $k \geq 0$ and any language $\LL$ we denote
by $\LL\chop_{k}$ (``L chop k'') the language obtained from $\LL$ by
performing, in the given order:

%%%%%%%%%%

\begin{enumerate}
\item
replace every $w \in \LL$ by $w\chop_{k}\,$, and
\item
add all subwords (= factors) to the language.
\end{enumerate}

The following properties of (laminary) languages are rather useful;
they follow directly from the definition.

\begin{lem}
\label{chopoff}
(a) Every laminary language $\LL$ satisfies, for every integer $k \geq
0$, the equality $\LL = \LL\chop_{k}$.

\smallskip
\noindent
(b) For every infinite language $\LL$ and for every integer $k$,
$L(\LL\chop_k)=L(\LL)$ and $\LL(L(\LL))=\cap_{k\in\N} \,\LL\chop_k$.
\qed
\end{lem}

Recall that a symbolic lamination $L \in \Lambda_{\CA}$ is {\it
minimal} if $L$ is equal to the closure of any of its orbits, with
respect to both, shift and inversion. This is equivalent to saying
that $L$ does not contain a proper sublamination.  One can easily
characterize laminary languages of such a minimal lamination:

\begin{defn}
A language $\cal L$ has the {\em bounded gap property} if for any word
$u$ in $\cal L$ there exists an integer $n=n(u) \in \N$ such that any
word $w \in \LL$ of length greater than $n$ contains $u$ or $u\inv$ as
a subword.
\end{defn}

The following is part of symbolic dynamics folklore \cite{fogg}:

\begin{prop}
\label{rem:boundedgap}
A (symbolic) lamination is minimal if and only if its laminary
language has the bounded gap property.  \qed
\end{prop}

Note that, if in addition the lamination is non-orientable, then for $n$
big enough any word $w$ of the laminary language will contain both, $u$
and $u\inv$.

%%%%%%%%%%%%%%%%%%%%%%%%%%%%%%%%%%%%%%%%%%%%%%%%%%%%%%%%%%%%%%%%

\section{Metrics and topology on the set of laminations}
\label{subsec:topologyonLambda}

For any laminary languages $\LL, \LL' \in \Lambda_{\LL}$ we define:
\[
d({\cal L},{\cal L}')=\exp(-\max(\{n\geq 0\ |\ {\cal L}_{ 2n+1} ={\cal
L}'_{2n+1}\} \cup \{0\})).
\]
This defines a distance on $\Lambda_{\LL}$ which is easily seen to be
ultra-metric, and it is clear that $\Lambda_\LL$ is a compact
Haussdorf totally disconnected perfect metric space: a Cantor set.

\smallskip

Similarly, one can define on the set $\Sigma_{\CA}$ of biinfinite
reduced words in $\CA^{\pm 1}$ a metric, by defining for any $Z, Z'
\in \Sigma_{\CA}$ the distance
\[
d(Z, Z')=\exp(-\max(\{n\geq 0\ |\ Z_n=Z'_n\}\cup \{0\})) \, ,
\]
where for any reduced biinfinite word $Z = \ldots z_{i-1} z_{i}
z_{i+1} \ldots$ we denote the {\em central subword of length $2n + 1$}
by $Z_{n} = z_{-n} z_{-n + 1} \ldots z_{n} \, $.

\smallskip

From these definitions and the shift-invariance of a symbolic
lamination we obtain directly that a symbolic lamination $L_{\CA}$ is
contained in the $\epsilon$-neigh\-bor\-hood in $\Sigma_{\CA}$ of a second
symbolic lamination $L'_{\CA}$ if and only if $\LL_{2n+1}(L_{\CA})$ is
a subset of $\LL_{2n+1}(L'_{\CA})$, for $\epsilon = e^{-n}$.  This
metric on $\Sigma_{\CA}$ induces a Hausdorff metric on the set
$\Lambda_{\CA}$ of symbolic laminations in $\CA^{\pm 1}$.  We obtain
directly:

\begin{prop}
\label{prop:isometrylanglamsymb}
The bijection $\rho^\CA_{\LL}: \Lambda_{\CA} \to \Lambda_{\LL}$ given
by $L_{\CA} \mapsto \LL(L_{\CA})$ is an isometry with respect to the
above defined metrics:
\[
d(L_{\CA}, L'_{\CA}) \leq e^{-n} \iff \LL_{2n+1}(L_{\CA}) =
\LL_{2n+1}(L'_{\CA})
\]
\qed
\end{prop}

As indicated in \S\ref{sec:alglam}, the choice of a basis $\CA$ of
the free group $\FN$ defines a word metric on $\FN$ and also a
(ultra-){\em metric at infinity} on $\partial \FN$, by specifying for
any $X, Y \in \partial \FN$, with prefixes $X_{n}$ and $Y_{n}$
respectively, the distance
\[
d_{\CA}(X,Y)=\exp(-\max\{n\geq 0\ |\ X_n=Y_n\}).
\]

In a similar vein as above for $\Sigma_{\CA}$, this distance can be
used to define a distance on $\partial^{2} \FN$, and we can define a
Hausdorff metric $d_{\CA}$ on $\Lambda^{2}(\FN)$. With a little care
we can show that this makes the bijection $\rho^{2}_{\CA}
:\Lambda^2(\FN) \to \Lambda_{\CA}$ from Proposition
\ref{prop:laminationflow} an isometry.  However, contrary to the case
of $\Lambda_{\CA}$ and $\Lambda_{\LL}$, the choice of a basis in $\FN$
and hence of the metric on $\partial \FN$ is not really natural, so
that we prefer for $\Lambda^2(\FN)$ only to consider the topology
induced by these metrics. Whenever a basis is specified, it is in any
case more convenient to pass directly to $\Lambda_{\CA}$ or to
$\Lambda_{\LL}$.  It is well known (and can easily be derived from the
material presented in \S\ref{sec:ccb} below) that different bases of
$\FN$ induce H\"older-equivalent metrics on $\partial \FN$ and on
$\partial^{2}\FN$, and thus also on $\Lambda^2(\FN)$. Thus we obtain:

\begin{prop}
\label{prop:homeos}
The canonical bijections
\[
\Lambda^2(\FN) \overset{\rho_{\CA}^{2}}{\longrightarrow}
\Lambda_{\CA} \overset{\rho^{\CA}_\LL}
{\longrightarrow} \Lambda_{\LL}
\]
are homeomorphisms.  They also preserve the partial order structure
defined on each of them by the inclusion as subsets.  \qed
\end{prop}

The topology on the space of laminations is explicitly encapsulated in
the following:

\begin{rem}
\label{convergence}
A sequence $(L^{2}_k)_{k\in\N}$ of algebraic laminations converges to
an algebraic lamination $L^{2}$ if and only if, for some (and hence
any) basis $\cal A$ of $\FN$, the sequence of corresponding symbolic
laminations $L_{k} = \rho_{\CA}^{2}(L^{2}_{k})$ and their presumed
limit $L = \rho_{\CA}^{2}(L^{2})$ satisfy the following:

\smallskip
\noindent
{\em Convergence criterion:} For any integer $n \geq 1$ there exists a
constant $K(n) \geq 1$ such that for all $k\geq K(n)$ one has:
\[
{\cal L}_n(L_{k})={\cal L}_n(L)\, .
\]
\end{rem}

\medskip

The following lemma will be used in \cite{chl1-III}.

\begin{lem}
\label{lem:sublaminations}
For any given algebraic lamination $L^{2}$ the set $\delta(L^{2})$ of
sublaminations of $L^{2}$ is a compact subset of $\Lambda^{2}$.
\end{lem}

\begin{proof}
Since $\Lambda^{2}$ is compact, it suffices to show that
$\delta(L^{2})$ is closed.  Any sublamination of $L^{2}$ has as
laminary language a sublanguage of the laminary language $\LL(L^{2})$
defined by $L^{2}$, and conversely.  Moreover, for laminary languages
the analogous statement as given by the lemma is trivially true, as
follows directly from the above Convergence criterion.
\end{proof}

We would like to point the reader's attention to the fact that the
space $\Lambda^{2}$ is rather large, and for some purposes perhaps too
large: it contains more objects than one would naturally think of as
analogues of surface laminations.  Of particular interest seems to be
the natural subspace of $\Lambda^{2}$ given by the closure $\overline
\Lambda_{\mbox{\scriptsize rat}} = \overline
\Lambda_{\mbox{\scriptsize rat}}(\FN)$ of the the space
$\Lambda_{\mbox{\scriptsize rat}}$ of rational laminations (compare
\S\ref{sec:alglam}).  We can now restate and prove Proposition
\ref{prop:notdense}:

\begin{prop}
\label{notdense2}
The inclusion $\overline \Lambda_{\mbox{\scriptsize rat}}
\subset \Lambda^2(\FN)$, for $N \geq 2$, is not an equality.
\end{prop}

\begin{proof}
For $a$ and $b$ in $\CA$ consider the symbolic lamination $L(\LL(Z))$
generated by the biinfinite word $Z = \ldots aaa \cdot bbb\ldots$. It
consists precisely of the $\sigma$-orbit of $Z$ and of the two
periodic words $\ldots aaa \cdot aaa\ldots$ and $\ldots bbb\cdot
bbb\ldots$, together with all of their inverses.  The laminary
language $\LL_{n}(Z)$ consists of the words $a^n, a^{n-1}b,
a^{n-2}b^{2}, \ldots, ab^{n-1}, b^n$ and their inverses. However,
every rational lamination $L$, with the property that the
corresponding laminary language contains these words, must contain the
rational sublamination $L(w)$ for some $w \in F(a, b)$ that contains
both letters, $a$ and $b$, or their inverses.  But then $\LL_{n}(L)$
must also contain the word $bx$ in $\LL_{2}(L)$, for some $x \in \CA
\cup \CA^{-1} \smallsetminus \{b, b\inv\}$.  This contradicts the
above Convergence criterion from Remark \ref{convergence}, for any
$L_{k} = L$ as above.
\end{proof}

On the other hand, the closure of the rational laminations seems to be
a reasonable subspace of $\Lambda^2$, as shown by the following:

\begin{prop}
\label{limitset}
$\overline \Lambda_{\mbox{\scriptsize rat}}$ contains all minimal
algebraic laminations.
\end{prop}

\begin{proof}
We prove the proposition for non-orientable minimal laminations, where
$\FN$-orbits and $< \FN, \hbox{\rm flip}>$-orbits agree, and leave the
generalization for orientable laminations to the reader.
 
Let $L^2$ be a minimal algebraic lamination and $\cal A$ a basis of
$F$. Let $L_{\CA} = \rho^{2}_{\CA}(L^{2})$ be the symbolic lamination
and ${\cal L}= \rho^\CA_{\LL}(L_{\cal A})$ the laminary language
canonically associated to $L^{2}$. By minimality of $L^{2}$ the
language ${\cal L}$ has the bounded gap property (see Proposition
\ref{rem:boundedgap}): For any integer $n$ there exists a bound
$K=K(n)$ such that for any words $u$ and $w$ of ${\cal L}$ where the
length of $u$ is smaller than $n$ and the length of $w$ is greater
than $K$, $u$ occurs as a subword of $w$.

This proves that for any word $w$ of $\cal L$ of length greater than
$K$ we have ${\cal L}_n(w)={\cal L}_n(L^2)$.  If moreover $w$ is
cyclically reduced, we obtain:
\[
{\cal L}_n(L(w))\supset{\cal L}_n(w)={\cal L}_n(L^2)
\]

Now let $u$ be any word of $\cal L$ of length $n$ and $v$ another word
of $\cal L$ of length $3K$. Write $v=w_1w_2w_3$ where $w_1$, $w_2$,
$w_3$ are all of length $K$: The product $w_1w_2w_3$ is reduced, and
each $w_{i}$ is a subword of $v$.  Now $u$ must be a subword of both,
$w_1$ and $w_3$: We can write the corresponding reduced products
$w_1=w'_1uw''_1$ and $w_3=w'_3uw''_3$, and we define:
\[
v' = u w''_1 w_2 w'_3
\]
Since $v'$ contains $w_{2}$ as subword, its length is bigger than $K$,
and hence the previous equality applies: ${\cal L}_n(v')={\cal
L}_n(L^{2})$. Moreover, since $w'_{3} u$ is a subword of the reduced
word $w_{3}$, it follows that $v'$ is cyclically reduced, and hence
${\cal L}_n(L(v')) \supset {\cal L}_n(L^{2})$.  Finally, since $u$ has
length $n$, any subword of length $n$ of the reduced biinfinite word
$\ldots v' v' \cdot v' v' \ldots$ that is not a subword of $v'$ is
necessarily a subword of $w_{2} w'_{3} u$, and hence of $v$.  Hence we
get ${\cal L}_n(L(v')) \subset {\cal L}_n(L^{2})$ and thus
\[
{\cal L}_n(L(v')) = {\cal L}_n(L^{2})\, .
\]
Thus, for any integer $n$ we found a word $v' = v'(n) \in F(\CA)$ such
that the rational lamination $L(v'(n))$ satisfies ${\cal
L}_n(L(v'(n)))={\cal L}_n(L^2)$. Hence the Convergence criterion of
Remark~\ref{convergence} gives directly that $L(v'(n)) \overset{n \to
\infty}{\longrightarrow} L^{2}$.  ${}_{}$
\end{proof}

The two previous propositions imply directly Theorem
\ref{thm:ratdense}.

%%%%%%%%%%%%%%%%%%%%%%%%%%%%%%%%%%%%%%%%%%%%%%%

\section{Bounded cancellation}
\label{sec:ccb}

An important tool when dealing with more than one basis in a free
group $\FN$ is {\em Cooper's cancellation bound} \cite{coop}.  We
denote by $\vbar w \rvbar_{\CA}$ the length of the element $w \in \FN$
when written as reduced word in a basis $\CA$ of $\FN$.

\begin{lem}\label{lem:ccb}
\label{lem:bbtcooper}
Let $\alpha$ be an automorphism of a free group $\FN$ and let $\cal A$
be a basis of $\FN$.  Then there exists a constant $C \geq 0$ such
that, for any elements $u,v \in \FN$ with
\[
\vbar u \rvbar_{\CA} +
\vbar v \rvbar_{\CA} \, \, = \, \,
\vbar uv \rvbar_{\CA}
\]
(i.e. there is no cancellation in the product $uv$ of the reduced
words $u$ and $v$)
one has
\[
0 \leq \vbar \alpha(u) \rvbar_{\CA} + \vbar \alpha(v) \rvbar_{\CA} -
\vbar \alpha(uv) \rvbar_{\CA} \, \, \leq \, \, 2 C
\]
\end{lem}

As any second base $\cal B$ is the preimage of $\CA$ under some
$\alpha \in \Aut(\FN)$, the last line of the above statement can
equivalently be replaced by
\[
0 \leq \vbar u \rvbar_{\cal B} + \vbar v \rvbar_{\cal B} - \vbar uv
\rvbar_{\cal B} \, \, \leq \, \, 2 C
\]
We denote by $\BBT(\CA,\alpha)$ or $\BBT(\CA,{\cal B})$ the smallest
such constant $C$.

\smallskip

An elementary proof of the above lemma can be given inductively, by
decomposing the given automorphism (or basis change) into elementary
Nielsen transformations.  In modern geometric group theory language,
one can restate the lemma as a special case of the fact that any two
word metrics on a group $G$ based on two different finite generating
systems give rise to a quasi-isometry which realizes the identity on
$G$.

\medskip

This lemma has been interpreted and generalized in term of maps
between trees in \cite{gjll}. We describe now this interpretation; a
generalization is given in \cite{chl1-II}.

Let $T_{\cal A}$ and $T_{\cal B}$ be the metric realisations (with
constant edge length 1) of the Cayley graphs of $\FN$ with respect to
$\cal A$ and $\cal B$.  Let $i = i_{\cal A, B}$ the equivariant map
from $T_{\cal A}$ to $T_{\cal B}$ which is the identity on vertices
and which is linear (and thus locally injective) on edges.  Then
Cooper's cancellation lemma~\ref{lem:bbtcooper} can be rephrased as:

\begin{lem}
For any (possibly infinite) geodesic $[P,Q]$ in $T_{\cal A}$ the image
$i([P,Q])$ lies in the $C$-neighborhood in $T_{\CB}$ of $[i(P),i(Q)]$,
for some $C > 0$ (in particular for $C = \BBT(\CA,{\cal B})$ as above)
independent on the choice of $P, Q \in T_{\CA}$.  \qed
\end{lem}

Finally, we state the following lemma that is used in \cite{chl1-II}:

\begin{lem}\label{lem:BBTconj}
Let $\A$ and $\cal B$ be two bases of $\FN$.  Any element $w$ of $\FN$
which is cyclically reduced with respect to the basis $\A$ is ``almost
cyclically reduced with respect to $\cal B$''.  More specifically, if
\[
w=y_1\cdots y_r y_{r+1}\cdots y_n y_r\inv\cdots y_1\inv
\]
with $y_{i} \in \CB^{\pm 1}$ is a reduced word (in particular with
$y_{r+1} \neq y_{r}\inv$ and $y_{n} \neq y_{r}$), then one has $r \leq
\BBT(\A,{\cal B})$.
\end{lem}

\begin{proof} 
Apply Lemma \ref{lem:bbtcooper} to $w^2$.
\end{proof}

%%%%%%%%%%%%%%%%%%%%%%%%%%%%%%%%%%%%%%%%%%%%%%%%%%%%%%%%%%%%%

\section{The $\Out(\FN)$-action on laminations and
laminary languages}
\label{sec:outfnactiononLambda}

In \S\ref{sec:alglam} we briefly mentioned that there is a natural
action by any automorphism of $\FN$ as homeomorphism on the boundary
$\partial\FN$, and thus on $\Lambda^2$.  This is a well known result
in geometric group theory: Indeed the very fact that the boundary of a
free group can be defined without any reference to a given basis is
exactly equivalent to that statement.  The key fact here is that a
basis change in $\FN$ (or, equivalently, an automorphism of $\FN$)
induces a change of the metric on $\FN$ (see
\S\ref{subsec:topologyonLambda}) in a Lipschitz equivalent way.
Therefore it changes the induced metric on the boundary (viewed as the
set of one-sided infinite reduced words, see
\S\ref{subsec:topologyonLambda}) in a H\"older equivalent way.

\medskip

A more direct combinatorial way to define the action of $\Out(\FN)$ on
languages is given as follows: Notice first that the elementwise image
$\alpha({\cal L})$ of a laminary language $\cal L$ under an
automorphism $\alpha \in \Aut (\FN)$ is in general not a laminary
language.

By Lemma~\ref{lem:ccb}, for $C=\BBT(\CA,\alpha)$ the language
$\alpha(\LL)\chop_C$ is laminary, and by Lemma \ref{chopoff} we have
$L(\alpha(\LL))=L(\alpha(\LL)\chop_C)$.  Thus, if we consider the
outer automorphism $\widehat \alpha \in \Out(\FN)$ defined by
$\alpha$, we can define:
\[
\widehat \alpha ({\cal L}) = \alpha(\LL)\chop_C={\cal
L}(L(\alpha({\cal L})))
\]
It follows directly from the second equality that this does not depend
on the choice of the automorphism $\alpha$ in the class
$\widehat\alpha$.  It also follows directly from our definitions that
this action of $\widehat \alpha$ is in fact a homeomorphism of the
space $\Lambda_{\LL}$ of laminary languages in $\CA^{\pm 1}$.

Similarly, for any symbolic lamination $L_{\cal A}$ we define
\[
\widehat \alpha (L_{\cal A}) = L(\alpha({\cal L}(L))) \, .
\]

From these definitions we see directly that the actions of
$\widehat\alpha$ commute with the (bijective) map
$\rho^{\CA}_{\LL}:\Lambda_{\CA} \to \Lambda_{\LL}$ given in
Proposition~\ref{bijectionsymbolic}.

If $\beta$ is a second automorphism of $\FN$ and $C'=
\BBT(\CA,\beta)$, one gets from Lemma \ref{lem:ccb} that
\[
\alpha(\beta(\LL)\chop_{C'})\chop_C=
(\alpha\beta)(\LL)\chop_{C''} \, ,
\]
with $C'' = \vbar \alpha \rvbar_{\CA}\, C' + C$ and $\vbar \alpha
\rvbar_{\CA} = \max\{ \vbar \alpha(x) \rvbar_{\CA}~: \, x \in \CA
\,\}$.  This shows that the definitions above give an action of
$\Out(\FN)$ on $\Lambda_\LL$ and on $\Lambda_\CA$.

\smallskip

Applying Lemma \ref{lem:ccb} again, we get that, if $(X,X')$ is a leaf
of an algebraic lamination $L^2$, then any subword of
$\rho_\CA(\alpha(X),\alpha(X'))$ is a word in $\alpha(\LL(X\inv
X'))\chop_C$.  This proves that $\rho^{2}_\CA$ is
$\Out(\FN)$-equivariant and thus concludes the proof of Theorem
\ref{theoremone}.

\medskip

Each of the above two versions of the $\Out(\FN)$-actions has its own
virtues: Surprisingly, the action on laminary languages generalizes
much more directly to more general homomorphisms $\phi: \FN \to F_{M}$
of free groups.  It is noteworthy in this context that non-injective
substitutions on biinfinite sequences are treated classically in
symbolic dynamics in a similar vein as injective ones, while from a
geometric group theory standpoint it is impossible to extend a
non-injective map $\phi$ as above in any meaningful way to a map
$\partial \phi: \partial \FN \to \partial F_{M}$.  The more common
injective case, however, is easy to understand even from the geometric
group theory standpoint:

\begin{rem}
\label{monomorphisms}
It is well known that every finitely generated subgroup of a free
group is quasi-convex. Thus an embedding $\phi: F_{M} \subset F_{N}$
induces canonically an embedding $\partial \phi: \partial \FM \subset
\partial \FN$, see \cite{gh}.  Clearly, this extends to an embedding
$\partial \phi^{2}: \partial^{2}\FM \subset \partial^{2}\FN$, but
since the image $\partial \phi^{2}(\partial^{2}\FM) \subset
\partial^{2}\FN$ is in general not $\FN$-invariant, an algebraic
lamination $L^{2} \subset \partial^{2}\FM$ is mapped by $\partial
\phi^{2}$ to a set $\partial \phi^{2}(L^{2}) \subset \partial^{2}\FN$
that is in general {\em not} an algebraic lamination. By taking the
closure of $\partial \phi^{2}(L^{2})$ with respect to the topology,
the $\FN$-action, and the flip map, one obtains however a well defined
algebraic lamination, which we denote by $\phi_{\Lambda}(L^{2})$, thus
defining a natural map:
\[
\phi_{\Lambda}: \Lambda^{2}(F_{M}) \to \Lambda^{2}(\FN)
\]
However, it has to be noted immediately that this map $\phi_{\Lambda}$
does not have to be injective: It suffices that the embedding $\phi$
maps elements $v,w \in \FM$ which are not conjugate in $\FM$ to
elements $\phi(v), \phi(w)$ that are conjugate in $\FN$: Then the
associated rational laminations satisfy
\[
L^{2}(v) \neq L^{2}(w) \in \Lambda^{2}(\FM) \, ,
\]
but also
\[
\phi_{\Lambda}(L^{2}(v))  = L^{2}(\phi(v)) = L^{2}(\phi(w))
= \phi_{\Lambda}(L^{2}(w))  \in \Lambda^{2}(\FN) \, .
\]

On the other hand, we note that if $\FM$ is a free factor of $\FN$, 
then the lamination space $\Lambda^{2}(\FM)$ is canonically embedded 
into $\Lambda^{2}(\FN)$: it suffices to consider a basis of $\FN$ 
which contains as a subset a basis of $\FM$.

It seems to be an interesting question of when precisely the map
$\phi_{\Lambda}: \Lambda^{2}(F_{M}) \to \Lambda^{2}(\FN)$
induced by an embedding $\phi: F_{M} \subset F_{N}$ is injective, 
and in particular, if this is the case if and only if the subgroup 
$\FM$ is malnormal in $\FN$.
\end{rem}

\bigskip

We finish this paper with an answer to the question we posed in
\S\ref{sec:alglam}.

\begin{prop}
\label{thierry}
Let $\CA$ be a basis of $\FN$, and let $a$ be an element of $\CA$.
Then, for any $N \geq 2$, the closure of the $\Out(\FN)$-orbit of the
rational lamination $L(a)$ in $\Lambda^{2}$ is not the only non-empty
minimal closed $\Out(\FN)$-invariant subspace of $\Lambda^2$.
\end{prop}

\begin{proof} 
Let $a$ be as above, and let $b$ be another element of $\CA$.
Consider the rational lamination $L([a,b])$. Then for any outer
automorphism $\hat\alpha$ of $\FN$ and any automorphism $\alpha$
representing it, one has
\[
\hat\alpha(L([a,b]))=L(\alpha([a,b])).
\]
As the derived subgroup is characteristic, the $\Out(\FN)$-orbit of
$L([a,b])$ consists of some minimal rational laminations associated to
cyclically reduced words of the derived subgroup. Now any cyclically
reduced word of the derived subgroup contains a subword of the form
$xy$, where $x,y$ are distinct elements of $\CA^{\pm 1}$ with $x \neq
y^{-1}$.  This proves that for any outer automorphism $\hat\alpha$,
the laminary language $\LL(\hat\alpha(L([a,b])))$ contains a reduced
word of the form $xy$.  It follows from the Convergence criterion in
Remark~\ref{convergence} that $L(a)$ is not in the closure of the
$\Out(\FN)$-orbit of $L([a,b])$.
\end{proof}

%\newpage
%%%%%%%%%%%%%%%%%%%%%%%%%%%%%%%%%%%

\affiliationone{
Thierry Coulbois, Arnaud Hilion and Martin Lustig\\
Math\'ematiques (LATP)\\
Universit\'e Paul C\'ezanne -- Aix-Marseille III\\
av. escadrille Normandie-Ni\'emen\\
13397 Marseille 20\\ 
France
\email{Thierry.Coulbois@univ-cezanne.fr\\
Arnaud.Hilion@univ-cezanne.fr\\
Martin.Lustig@univ-cezanne.fr\\
}
}

\end{document}